\documentclass[12pt]{article}

\usepackage{latexsym}
\usepackage{amsmath}
\usepackage{amsfonts}
\usepackage{amssymb}
\usepackage{amsthm}
\usepackage{amscd}
\usepackage{color}
\input xy
\xyoption{all} \hfuzz=1.5pt \tolerance=500

\newcommand{\ot}{\otimes}

\newcommand{\tr}{\triangleright}
\newcommand{\tl}{\triangleleft}

\newcommand{\va}{\varphi}

\newcommand{\om}{\omega}

\newcommand{\id}{{\bf 1}}
\newcommand{\co}{{\B C}}

\newcommand{\bd}{\begin{document}}
\newcommand{\ed}{\end{document}}
\newcommand{\Htp}{\mathbin{\stackrel{\cdot}{\bigotimes}}}
\newcommand{\Hts}{\mathbin{\stackrel{\cdot}{\bigoplus}}}
\newcommand{\nl}{\nu\in\Lambda}
\newcommand{\xl}{X_\nu;\nl}
\newcommand{\x}{X_\nu}
\newcommand{\la}{\langle}
\newcommand{\ra}{\rangle}
\newcommand{\ch}{\cal H}
\newcommand{\su}{\subseteq}
\newcommand{\B}{\Bbb}
\newcommand{\e}{\varepsilon}
\newcommand{\ov}{\overline}
\newcommand{\vk}{\varkappa}
\newcommand{\vt}{\vartheta}
\newcommand{\al}{\alpha}
\newcommand{\de}{\delta}
\newcommand{\lo}{\Longleftrightarrow}
\newcommand{\cK}{{\cal F}}
\newcommand{\sq}{{\square}}
\newcommand{\bb}{{\cal B}}
\newcommand{\dd}{{\cal D}}
\newcommand{\rr}{{\cal R}}
\newcommand{\cR}{{\cal R}}
\newcommand{\f}{{\cal K}}
\newcommand{\cb}{{\cal CB}}
\newcommand{\wrr}{\widetilde{\cal R}}
\newcommand{\br}{^\bullet{\cal R}}
\newcommand{\en}{E_\nu}
\newcommand{\el}{{\cal L}}
\newcommand{\lr}{\Longrightarrow}
\newcommand{\bp}{\bigcirc}
\newcommand{\Long}{\Longleftarrow}
\newcommand{\q}{\quad}
\newcommand{\qq}{\qquad}
\newcommand{\cd}{\cdot}
\newcommand{\fur}{f: E\to{\B R}}
\newcommand{\furo}{f_0: E_0\to{\B R}}
\newcommand{\fuc}{f: E\to{\B C}}
\newcommand{\ii}{\infty}
\newcommand{\di}{\diamondsuit}
\newcommand{\bgd}{{\bigtriangledown}}
\newcommand{\bu}{{\bigtriangleup}}
\newcommand{\bc}{{completely bounded}}
\newcommand{\cc}{{completely contractive}}
\newcommand{\qs}{{quantum space}}
\newcommand{\isc}{{isometric}}
\newcommand{\ism}{{isomorphism}}
\newcommand{\qss}{{quantum spaces}}
\newcommand{\bco}{{completely bounded operator}}
\newcommand{\bcos}{{completely bounded operators}}
\newcommand{\res}{{respectively}}
\newcommand{\tp}{{tensor product}}
\newcommand{\eq}{{equivalent}}
\newcommand{\qtp}{{quantum tensor product}}
\newcommand{\mma}{\mathrel{\mathop{\otimes}\limits_{A}}}
\newcommand{\mmb}{\mathrel{\mathop{\otimes}\limits_{\bb}}}
\newcommand{\mmp}{\mathrel{\mathop{\otimes}\limits_{p}}}
\newcommand{\mmh}{\mathrel{\mathop{\otimes}\limits_{h}}}
\newcommand{\mmop}{\mathrel{\mathop{\otimes}\limits_{4}}}
\newcommand{\mmf}{\mathrel{\mathop{\otimes}\limits_{op}}}
\newcommand{\mmi}{\mathrel{\mathop{\otimes}\limits_{i}}}
\newcommand{\mmm}{\mathrel{\mathop{\otimes}\limits_{\bb-\bb}}}
\newcommand{\mmdd}{\mathrel{\mathop{\otimes}\limits_{\cdot}}}

\newcommand{\mmo}{\mathrel{\mathop{\cdot}\limits_{1}}}
\newcommand{\mmA}{\mathrel{\mathop{\otimes}\limits_{A-A}}}
\newcommand{\mmd}{\mathrel{\mathop{\cdot}\limits_{2}}}
\newcommand{\msp}{\mathrel{\mathop{\otimes}\limits_{sp}}}
\newcommand{\mms}{\stackrel{h}{\otimes}}
\newcommand{\mmt}{\stackrel{4}{\otimes}}
\newcommand{\mmx}{\stackrel{i}{\otimes}}
\newcommand{\mmy}{\stackrel{p}{\otimes}}
\newcommand{\mmz}{\stackrel{sp}{\otimes}}
\newcommand{\hil}{\stackrel{\cd}{\otimes}}
\newcommand{\gd}{\ddagger}
\newcommand{\od}{\odot}
\newcommand{\mt}{\mapsto}
\newcommand{\mmc}{\mathrel{\mathop{\otimes}\limits_{\sim}}}
\newcommand{\mme}{\stackrel{\sim}{\otimes}}

\begin{document}


\bigskip

\centerline{{\bf\large Projective quantum modules and projective ideals of $C^*$-algebras
}}

   \vspace{1cm}

\centerline{A.~Ya.~Helemskii}

\centerline{Faculty of Mechanics and Mathematics}

\centerline{Moscow State (Lomonosov) University}

\centerline{helemskii@rambler.ru}

\vspace{1cm}

{\bf Abstract}

\bigskip
{\small We introduce in non-coordinate presentation the notions of a quantum algebra and of a quantum
 module over such an algebra. Then we give the definition of a projective quantum module and of a free
 quantum module,
 the latter as a particular case of the notion of a free object in a rigged category.
 (Here we say ``quantum'' instead of frequently used protean adjective ``operator''). After this we discuss
 the general
 connection between projectivity  and freeness. Then we show that for a Banach quantum algebra
 $A$ and a Banach quantum space  $E$ the Banach quantum $A$-module
 $A\widehat\ot_{op}E$ is free, where `` $\widehat\ot_{op}$ '' denotes the operator-projective
tensor product of Banach quantum spaces. This is used in the proof of the following theorem: all closed
left ideals in a separable $C^*$-algebra, endowed with the standard quantization, are projective left quantum
 modules over this algebra.

 \bigskip
 Bibliography: 29 titles. }

\vspace{1cm}


\vspace{1cm}

{\bf Introduction}

\bigskip
As it is well known, the concept of projective module is very important in algebra, and it is one of the
three pillars, on which the whole building of homological algebra rests. (Two others are the notions of an
injective and a flat module). After the sufficient development of the theory of Banach and operator algebras
 the concept of the projectivity was carried over to this area, at first in the context of ``classical'' and
 later of ``quantum'' functional analysis; see, e.g.,~\cite{he1,heb1,wood})
  By the latter we mean the area more frequently called (abstract) operator space theory.

The present paper consists of three parts. In the first part, after necessary preparations, we discuss
the notion of a projective quantum (= operator) module over a quantum (= operator) algebra.
Actually, there are quite a few different approaches to what to call projective module, quantum as well as
``classical''. We concentrate on the so-called {\it relative} projectivity, which is most known and
developed in the ``classical'' context. As to other existing versions, we just mention two of them, called
topological and metric.
Our definitions are given in the frame-work of the so-called non-coordinate, and not of the
more widespread ``matrix'', approach to the notion of the
operator space. The latter approach is presented in the widely known textbooks~\cite{efr,pis,paul,blem}.  
 The non-coordinate presentation is, in our opinion, more convenient for this circle of questions, which is intimately connected with tensor products.

In the second part we discuss a method of verifying whether a given quantum module is projective. This
 method is based on the notion of the so-called freeness.
It has a general character and can be applied to the broad variety of versions of projectivity,
appearing in algebra, functional analysis and topology. The needed definitions, gathered in~\cite{he2},
generalize those given by MacLane
in his theory of relative Abelian categories~\cite{mcl}. (Note that typical categories of functional
analysis we have to work with are never Abelian, and often even not additive). We want to emphasize that
all results of general-categorical character we use are actually particular cases or direct corollaries of
results on adjoint functors contained
in the book of MacLane~\cite{mc2}; we only present these facts in the language, suitable for our aims.
In~\cite{he2,he3} this categorical approach was applied to several versions of projectivity for
``classical'' normed modules and for the so-called metric projectivity of quantum modules. Now we apply it
for the relative projectivity of quantum modules.

In the concluding part of our paper we apply the mentioned method to $C^*$-algebras, endowed with 
their standard quantization. We show that in the case of separable algebras all their closed left 
ideals, being considered as quantum modules, are relatively projective. Such a result has some 
historical background. Many years ago, in the context of ``classical'' Banach modules, relatively 
projective closed ideals of commutative $C^*$-algebras were characterised as those with paracompact 
Gelfand spectrum~\cite{he4}. In particular, all closed ideals of separable commutative 
$C^*$-algebras are relatively projective. Later Z.~Lykova~\cite{lyk} has proved that the latter 
assertion is valid for all closed left ideals in an arbitrary separable $C^*$-algebra; the similar 
theorem for an important particular case was obtained by J.Phillips/I.Raeburn~\cite{phr}. Therefore 
our result can be considered as a ``quantum'' version of the theorem of Lykova. 

\bigskip
{\bf I. Preliminaries}

\bigskip
The non-coordinate presentation of quantum functional analysis (= operator space theory) is the subject
 of the book~\cite{heb2}. Nevertheless, for the convenience of the reader, we shall briefly recall three
 most needed definitions.

\medskip
To begin with, let us choose an arbitrary separable infinite-dimensional Hilbert space, denote it
by $L$ and fix it throughout the whole paper. As usual, by $\bb(E,F)$ we denote the space of all
bounded operators between respective normed spaces with the operator norm. We write $\bb$ instead
of $\bb(L,L)$.

The symbol $\ot$ is used for the (algebraic) tensor product of linear spaces and, unless stated 
explicitly otherwise, for elementary tensors. The symbols $\ot_p$ and $\ot_i$ denote the 
non-completed projective and, respectively, injective tensor product of normed spaces whereas the
symbols $\widehat\ot_p$ and $\widehat\ot_i$ denote the respective completed versions of these 
tensor products (cf., e.g.,~\cite{rya}). The symbol $\hil$ is used for the Hilbert tensor product 
of Hilbert spaces and bounded operators, acting on these spaces. 

Denote by $\cK$ the (non-closed) two-sided ideal of $\bb$, consisting of finite rank bounded 
operators. Recall that, as a linear space, $\cK=L\ot L^{cc}$, where $(\cd)^{cc}$ is the symbol of 
the complex-conjugate space. More precisely, there is a linear isomorphism $L\ot L^{cc}\to\cK$, 
well-defined by taking $\xi\ot\eta$ to the operator $\zeta\mt\la\zeta,\eta\ra\xi$. 

The  operator norm on $\cK$ is denoted just by $\|\cd\|$; it corresponds, after the indicated 
identification, to the norm on $L\ot_i L^{cc}$. When we say ``the normed space $\cK$'', we mean, 
unless mentioned explicitly otherwise, the operator norm. However, sometimes we shall need also the 
trace-class norm $\|\cd\|_{\cal N}$ on $\cK$, corresponding to the norm on  $L\ot_p L^{cc}$. The 
space $\cK$ with that norm will be denoted by ${\cal N}$. 

If $E$ is a linear space, the identity operator on $E$ will be denoted by $\id_E$. We write just $\id$
instead of $\id_L$.

\medskip
The basic concepts of the theory of operator spaces are based on the triple notion of the 
amplification, first of linear spaces, then of linear operators and finally of bilinear operators. 

The {\it amplification} of a given linear space $E$
 is the tensor product $\cK\ot E$. Usually we shall briefly denote it by $\cK E$, and an elementary
tensor $\xi\ot x; \xi\in L, x\in E$, by $\xi x$.

\medskip
{\bf Remark.} In this way we behave according to the general philosophy of quantum or 
non-commutative mathematics. Indeed, we take a definition of a basic notion of an area in question 
and replace in it ``a thing commutative'' by ``a thing non-commutative''. In our case we replace 
complex scalars in the definition of a linear space $E=\co E$ by ``non-commutative scalars'' from 
$\cK$). 

\medskip
The important thing is that $\cK E$ is a bimodule over the algebra $\bb$ with the outer multiplications,
well defined by $a\cd(bx):=(ab)x$ and $(ax)\cd b:=(ab)x$. An (ortho)projection $P\in\bb$ is called
a {\it support} of an element $u\in\cK E$, if we have $P\cd u=u=u\cd P$.

\medskip
{\bf Definition 1.} A norm on $\cK E$ is called {\it abstract operator norm}, or, as we prefer to say,
{\it quantum norm} (for brevity, Q-norm) {\it on $E$}, if it satisfies two conditions, the so-called Ruan's
axioms:

\smallskip
(i) the $\bb$-bimodule $\cK E$ is contractive, that is we always have the estimate $\|a\cd u\cd
b\|\le\|a\|\|u\|\|b\|$; $a,b\in\bb, u\in\cK E$

\smallskip
(ii) if $u,v\in\cK E$ have orthogonal supports, then $\|u+v\|=\max\{\|u\|,\|v\|\}$

\smallskip
A space $E$, endowed by a Q-norm, is called an {\it abstract operator space}, or a {\it quantum 
space} (for brevity, a {\it Q-space}). 

\medskip
In an obvious way, every subspace of a Q-space also becomes a  Q-space.

Note that a Q-space $E$ becomes a usual (``classical'') normed space, if for $x\in E$ we set 
$\|x\|:=\|px\|$, where $p$ is an arbitrary rank 1 projection; by virtue of the axiom (i), the 
number $\|x\|$ does not depend on a particular choice of $p$ and obviously satisfies the definition 
of a norm. The resulting normed space is called the {\it underlying space} of a given Q-space 
whereas the latter is called a {\it quantization} of the former. 

As it is shown in~\cite[Prop. 2.2.4]{heb2}, we always have
$$
\|ax\|=\|a\|\|x\|,\eqno(1)
$$
in other words, a Q-norm on $E$ as a norm on $\cK\ot E$ is a cross-norm.

Obviously, the simplest space $\co$ has the only quantization, obtained by the identification of
$\cK\co$ with $\cK$.

We refer to the cited textbooks for numerous examples of Q-spaces, including the most important 
and, in a sense, universal recipe of the quantization of a space, consisting of operators. However, 
in the present paper, we only need 

\medskip
{\bf Example.} Let $A$ be a $C^*$-algebra. In this case $\cK A$, as a tensor product
of involutive algebras, is
itself an involutive algebra. Moreover, $\cK$ is obviously a union of an increasing net of
finite-dimensional, and hence nuclear, $C^*$-algebras. From this one can easily observe
(cf.~\cite[Section 2.3]{heb2}) that the algebra $\cK A$ has a unique norm, possessing the $C^*$-property,
and this (non-complete) norm is a Q-norm on $A$. In what follows, the latter Q-norm will be called
{\it standard}.

\bigskip
Now suppose that we are given an operator $\va:E\to F$ between linear spaces. Denote, for brevity,
by $\va_\ii$ the operator $\id_\cK\ot\va:\cK E\to\cK F$, well defined on elementary tensors by
$ax\mt a\va(x)$, and call it the {\it amplification} of $\va$. Clearly, $\va_\ii$ is a
morphism of $\bb$-bimodules.

\medskip
{\bf Definition 2.} An operator $\va:E\to F$ between Q-spaces is called {\it completely bounded},
respectively, {\it completely contractive}, {\it completely isometric, completely isometric
isomorphism},  if $\va_\ii$ is bounded, respectively contractive, isometric, isometric isomorphism.
We write $\|\va\|_{cb}:=\|\va_\ii\|$ and call it {\it completely bounded norm of} $\va$.

\medskip
If $\va$ is bounded in the ``classical'' sense, that is being considered between respective
underlying normed spaces, we say that it is (just) {\it bounded} and denote its operator norm,
as usual, by $\|\va\|$. It is easy to see that every completely bounded operator $\va:E\to F$ is
obviously bounded, and we have $\|\va\|\le\|\va\|_{cb}$ .

As to various examples of completely bounded operators, as well of bounded not completely bounded
operators see, e.g.,~\cite[Section 3.2]{heb2}). We only note that every involutive homomorphism
between $C^*$-algebras is ``automatically'' completely contractive~\cite[Theorem 3.2.10]{heb2}; this is
the non-coordinate presentation of what was said in~\cite{efr}.

As in the classical analysis, among Q-spaces those that are complete seem the most important. We
say that a normed Q-space is {\it complete} (or Banach), if its underlying normed space is
complete.

The  {\it completion} of a normed Q-space, say $E$, is by definition, is a pair 
$(\overline{E},i:E\to\overline{E})$, consisting of a complete Q-space and a completely isometric 
operator, and such that the same pair, considered for underlying spaces, is the ``classical''   
completion of $E$ as of a normed space. It is easy to see that for every normed Q-space $E$ there 
exists a completion. (The simple argument, is given in~\cite[Section 3]{heb2}). 
 Also it is easy to observe that the ``quantum'' completion has the universal property similar to that of
  the ``classical'' completion.
  Namely, if $(\overline{E}, i)$ is a
  completion of a Q-space $E$, $F$ a complete Q-space and $\va: E\to F$ a completely bounded operator,
  then there exists a unique completely bounded operator
  $\overline{\va}:\overline{E}\to F$, extending, in the obvious sense, $\va$. Moreover, we have
  $\|\overline{\va}\|_{cb}=\|{\va}\|_{cb}$.

\bigskip
Bilinear operators also can be amplified, however, in two essentially different ways. Namely, for a given
bilinear operator between linear spaces there are two standard ways to construct a bilinear operator
between respective amplifications. One of these constructions is called in~\cite[1.6]{heb2} strong
and another weak amplification. In the present paper we need only weak amplification, so we shall refer
to it as to (just) amplification.

To give the relevant definition, we need an operation  that would imitate the tensor multiplication of operators on our canonical Hilbert space $L$ but would not lead out of this space. For this aim, we supply $L$ by a sort of additional structure.

By virtue of Riesz-Fisher Theorem, there exists plenty of unitary isomorphisms between Hilbert 
spaces $L\hil L$ and $L$. Choose and fix one, say, $\iota:L\hil L\to L$, throughout our whole 
paper. After this, for given $\xi,\eta\in L$ we denote the vector \\ $\iota(\xi\ot\eta)\in L$ by 
$\xi\di\eta$, and for given $a,b\in\bb$ we denote the operator $ \iota(a\hil b)\iota^{-1}:L\to L$ 
by $a\di b$; obviously, the latter is well defined by the equality $(a\di b)(\xi\di\eta)= a(\xi)\di 
b(\eta)$. Also it is evident that we have the identities 
$$
(a\di b)(c\di d)=ac\di bd, \eqno (2)
$$
$$
\|\xi\di\eta\|=\|\xi\|\|\eta\| \qq{\rm and}\qq \|a\di b\|=\|a\|\|b\|. \eqno (3)
$$

\medskip
Now let $\rr:E\times F\to G$ be a bilinear operator between linear spaces. Its {\it amplification}
is the bilinear operator $\rr_\ii:\cK E\times\cK F\to\cK G$, well defined on elementary tensors by
$\rr_\ii(ax,by)=(a\di b)\rr(x,y)$.


\medskip
{\bf Definition 3}. Bilinear operator $\rr$ between Q-spaces is called {\it completely
bounded}, respectively, {\it completely contractive}, if its amplification is (just) bounded,
respectively, contractive. The norm of the latter amplification is called {\it completely bounded
norm} of $\rr$ and denoted by $\|\rr\|_{cb}$.

\medskip
As to numerous examples and counter-examples see, e.g.,~\cite[Section 5.2]{heb2} .

\medskip
For our future aims, we need one more version of the ``diamond multiplication''. Namely,
 for a linear space $E, a\in\cK$ and $u\in\cK E$ we introduce in $\cK E$ the elements, denoted by
 $a\di u$ and $u\di a$. They are well  defined by  assuming that the operation $\di$ is additive on both
  arguments and setting, for elementary tensors, $a\di bx:=(a\di b)x$ and $bx\di a:=(b\di a)x$.
 As it was shown in~\cite[Prop. 2.2.6]{heb2},
 we always have
$$
\|a\di u\|=\|a\|\|u\|=\|u\di a\|.\eqno(4)
$$

From now on we are already outside the scope of~\cite{heb2}.

In what follows, we shall often need some formulae, connecting some elements of amplifications of spaces
with some linear and/or bilinear operators. As a rule, these formulae can be easily verified on
elementary tensors and then, by additivity, extended to theirs sums, that is to general elements.
To avoid tiresome repetitions, in such cases we shall write just ``LOOK AT ELEMENTARY TENSORS''.

Recall that a bilinear operator $\rr:E\times F\to G$ gives rise, for every $x\in E$ and $y\in F$
to linear operator $^x\rr:F\to G:y\mt\rr(x,y)$ and $\rr^y:E\to G:x\mt\rr(x,y)$, sometimes called 
partial. For our future aims, let us notice 

\medskip
{\bf Proposition 1} (cf. also~\cite{ne1}). {\it If $E,F,G$ are Q-spaces and $\rr$ is completely bounded, then
for every
$x\in E$ and $y\in F$ the operators $^x\rr$ and $\rr^y$ are completely bounded. Moreover, we have
$\|^x\rr\|_{cb}\le\|x\|\|\rr\|_{cb}$ and $\|\rr^y\|_{cb}\le\|y\|\|\rr\|_{cb}$.}

\smallskip
$\tl$ Take an arbitrary $a\in\cK, x\in E, y\in F, u\in\cK E, v\in\cK F$ and observe the formula
$$
\rr_\ii(ax, v)=a\di[^x\rr_\ii(v)].
$$
(LOOK AT ELEMENTARY TENSORS). Combining the latter with (4) and (1), for every $a\in\cK;\|a\|=1$ we have
$$
\|^x\rr_\ii(v)\|=\|a\di[^x\rr_\ii(v)]\|=\|\rr_\ii(ax, v\|\le\|\rr\|_{cb}\|ax\|\|v\|=\|\rr\|_{cb}\|x\|\|v\|.
$$
This proves that $^x\rr$ is completely bounded, together with the first estimate. A similar
argument proves the rest. $\tr$

\bigskip
So far we discussed the quantization of spaces; now we turn to (complex associative) algebras and their
modules.
 As a matter of fact, there are two essentially different definitions of what could be called quantum
 algebra. One is based on the notion of a strong amplification of a bilinear operator (cf. above)
 It gives rise to the class of algebras, which is the subject of a deep and well developed theory with
 mighty
 theorems, concerning the operator realization of these algebras~\cite{brs}~\cite{ble2}. It is presented in
 the book of D.Blecher/C.Le Merdy~\cite{blem}. However, in this paper we choose a somewhat larger class, based on the notion of what we call here just amplification  of a bilinear operator.


\medskip
{\bf Definition 4}. Let $A$ be an algebra and simultaneously a Q-space. We say that $A$ is a
{\it Q-algebra}, if the respective bilinear operator of multiplication is completely contractive.

\medskip
Here is our main example.

\medskip
{\bf Proposition 2.} {\it A $C^*$-algebra with the standard Q-norm is a Q-algebra.}

\smallskip
The proof, given in~\cite[Theorem 5.1.3]{heb2}, actually uses the connection between the strong and the weak amplifications of bilinear operators that is outside of the scope of the present paper. Since we do not consider in this paper the strong amplification, we shall give, for the convenience of the reader, a straightforward proof.

$\tl$ Our task is to show that the bilinear operator ${\cal M}_\ii:\cK A\times\cK A\to\cK A$,
where ${\cal M}:A\times A\to A$ is the bilinear operator of multiplication in $A$, is contractive.

As we remember (see Example), $\cK A$ is an involutive algebra with a norm, possessing $C^*$-property and
hence the multiplicative inequality. Further, for all $u,v\in\cK A, a,b\in\cK$ we  have,  by virtue of (2),
 the formula
$$
{\cal M}_\ii(u\cd a,b\cd v)=(u\di b)(a\di v)\eqno(5)
$$
(LOOK AT ELEMENTARY TENSORS). Finally, presenting $u$ and $v$ as sums of elementary tensors,
we easily see that they have the same support, say $P$, of finite rank. Therefore, taking in (3) $a:=b:=P$
and using (5) and (4), we see that
$\|{\cal M}_\ii(u,v)\|=\|{\cal M}_\ii(u\cd P,P\cd v)\|
\le\|u\di P\|\|P\di v\|=\|u\|\|v\|$. $\tr$


\medskip
{\bf Definition 5}. Let $A$ be a Q-algebra, $X$ a left $A$-module in algebraic sense and
simultaneously a Q-space. We say that $X$ is a {\it left Q-module}, if the respective bilinear operator of
outer multiplication
is completely contractive.

\medskip
As an important class of examples, it is obvious that every left ideal in a Q-algebra $A$, considered with the Q-norm of a subspace and with the inner multiplication in the capacity of an outer multiplication, is a left Q-module over $A$.

If the underling space of a Q-algebra or of a Q-module is complete, we speak of a Banach Q-algebra 
or, respectively, a Banach Q-module. 

When we speak about a morphism between Q-modules, we always mean a morphism in algebraic sense
which is completely bounded as an operator.

Fix, for a time, a Q-algebra, say, $A$. Suppose we have a left Q-module over $A$, say $P$,
 which, for some reason, arouses our interest. We associate with this module the so-called {\it lifting data},  consisting of the following two things: a surjective $A$-module morphism $\tau:Y\to X$ between some other left Q-modules over $A$, and an arbitrary $A$-module morphism from $P$ into $X$.
The {\it lifting problem} is to find an $A$-module morphism $\psi$,
may be with some additional properties, making the diagram
\[
\xymatrix@R-10pt@C+15pt{
& Y \ar[d]^{\tau}\\
P \ar[ur]^{\psi} \ar[r]^{\va} & X }
\]
\noindent commutative. Such a $\psi$ is called a {\it lifting} of $\va$ across $\tau$.

We shall mainly concentrated on a certain version of projectivity, that is the oldest and most 
known in the ``classical'' context (see, e.g.,~\cite{he1}). However, we shall present it in the 
quantum context. 

Let us call a morphism $\tau:Y\to X$ between Q-modules {\it admissible} if it
has a right inverse completely bounded operator (generally speaking, not morphism), that is $\rho$ with
$\tau\rho=\id_X$.

\medskip
{\bf Definition 6}. A Q-module $P$ is called {\it relatively projective}, if for every admissible
morphism $\tau:Y\to X$
of Q-modules and an arbitrary morphism $\va:P\to X$ the respective lifting problem has a solution.

If in this definition we suppose that all participating modules are Banach Q-modules over a Banach
Q-algebra, we obtain the definition of a {\it relatively projective Banach Q-module.}

\medskip
(In many papers, concerning the well known ``classical'' counterpart of this definition for Banach modules,
introduced in~\cite{he1},
some people say, instead of ``relatively projective'', ``traditionally projective'', and some
just ``projective'').

\medskip
{\bf Remark.} One of advantages of the relative projectivity is that this property can be equivalently
 expressed in the language of derivations. We shall not give here details. We only mention, rather vaguely,
 that a left Q-module $P$ over $A$ is relatively projective if, and only if  every
 completely bounded derivation of $A$ with values in a certain class of quantum
bimodules over $A$, defined in terms of $P$, is inner.

\medskip
{\bf Remark.} It is worth noting that the quantization of modules can make non-projective modules projective and vice versa. For example, take  an infinite-dimensional Hilbert space $H$ and the algebra ${\cal B}(H)$ Then, as it is proved in~\cite{he5}, the
$A$-module $H\stackrel{\cd}{\otimes} H$, equipped with the action, well defined by $a\cd (x\ot
y):=a(x)\ot y$, is not projective in the classical sense but becomes projective after some natural
quantization. On the other hand, the same $H$, as it is known long ago, is classically projective
as a module over the algebra ${\cal N}(H)$ of trace-class operators on $H$, with  the action $a\cd x:=a(x)$.
However, O.Aristov, embarking
from some observations in~\cite{blm2}, suggested
such a quantization of $A$ and $H$ that we obtain a non-projective quantum ${\cal N}(H)$-module
(see~\cite{he6}).

\medskip
As to other existing types of the projectivity, we shall just give two definitions. First, let us call a
completely bounded operator $\tau:F\to E$ between Q-spaces {\it completely open}, respectively,
{\it completely strictly coisometric}, if its amplification $\tau_\ii:\cK F\to\cK E$ is (just)
open, respectively, strictly coisometric. (``Strictly coisometric'' means that our operator
 maps the closed unit ball of the domain space onto the closed unit ball of the range space).

\medskip
{\bf Definition 7}. A Q-module $P$ is called {\it topologically projective}, respectively, {\it
metrically projective}, if for every completely open, respectively, completely strictly coisometric
morphism $\tau:Y\to X$ and an arbitrary completely bounded morphism $\va:P\to X$ the relevant lifting problem has some solution $\psi$, respectively, a solution $\psi$ with the additional property
$\|\psi\|_{cb}=\|\va\|_{cb}$.

\medskip
Again, there is an obvious version of both notions for Banach Q-modules.

\medskip
{\bf Remark.} It is obvious that topological projectivity implies relative projectivity. Also it is 
known that metric projectivity implies topological projectivity. This follows from the result of 
S.Shteiner~\cite[Prop. 2.1.5]{sht}), obtained with the help of methods, based on the notion of a 
free module (see the next section). On the other hand, both converse statements are false 
(N.Nemesh, oral communication). One can see this, considering just Q-spaces (that is the case 
$A:=\co$), endowed with the so-called maximal quantization. (The latter is defined and discussed, 
e.g., in~\cite{efr} or in~\cite{heb2}). 

\bigskip
{\bf 2. Projectivity in rigged categories and freeness}

\bigskip
We proceed to a general-categorical method to prove or disprove the projectivity of a given
module. It is based on the notion of freeness (cf. Introduction).

Let ${\cal  K}$ be an arbitrary category. A {\it rig} of ${\cal  K}$ is a faithful (that is, not
gluing morphisms) covariant functor $\square:{\cal K}\to{\cal L}$, where ${\cal L}$ is another
category. A pair, consisting of a category and its rig, is called a {\it rigged category}. If a rig
is given, we shall call ${\cal  K}$ the main, and ${\cal L}$ the auxiliary category.

Fix, for a time, a rigged category, say $({\cal K},\square:{\cal K}\to{\cal L})$. We call a
morphism $\tau$ in ${\cal K}$ {\it admissible}, if $\square(\tau)$ is a retraction (that is, has a
right inverse morphism) in ${\cal L}$. After this, we call an object $P$ in ${\cal K}$ {\it
$\sq$-projective}, if, for every admissible morphism $\tau:Y\to X$ and an arbitrary morphism
$\va:P\to X$ in ${\cal K}$, there exists a lifting (now in the obvious general-categorical sense)
of $\va$ across $\tau$.

Let us denote the category of Banach Q-spaces and completely bounded operators by ${\bf QBan}$
and the category of left Banach Q-modules over a Banach Q-algebra $A$ and their (completely
bounded) morphisms by ${\bf QA-mod}$. (Here and thereafter, just to be definite, we consider the ``complete'' case; the ``non-complete'' case can be considered with the obvious modifications.)

Now one can immediately see that

{\it a Banach Q-module over a Banach Q-algebra is relatively projective if,
and only if it is $\sq$-projective with respect to the rig
$$
\sq:{\bf QA\mbox{-}mod}\to{\bf QBan},\eqno (6)
$$
\noindent where $\sq$ is the relevant forgetful functor.}

(We mean, of course, that $\sq$ forgets about the outer multiplication).

\medskip
{\bf Remark.} The topological and the metric projectivity also can be described in terms of suitable rigged categories. Here we only mention, rather vaguely, that in the ``topological''
case the respective functor $\sq$ forgets not only about the outer multiplication, but even about the
additive structure, and in the ``metric'' case it forgets about everything (a sort of ``complete amnesia''), so that our auxiliary category is just  the category of sets. See details, concerning
 the topological and metric projectivity, in~\cite{sht} and~\cite{he6}, respectively.

\medskip
We turn to the freeness.
Actually, we obtain its definition (which must be well known, perhaps under different names), if we shall scrutinize the ancient classical example of a free object, the free group.
Consider an arbitrary rig $\square:{\cal K}\to{\cal L}$ and an object $M$ in the auxiliary category ${\cal L}$.  An object ${\bf Fr}(M)$ in ${\cal K}$ is called {\it free} (or, to be precise, {\it
$\square$-free) object with the base $M$}, if,  for every $X\in{\cal K}$, there exists a bijection
$$
{\cal I}_{M,X}: {\bf h}_{{\cal L}}(M,\square X)\to{\bf h}_{{\cal K}}({\bf Fr}(M),X), \eqno (7)
$$
between the respective sets of morphisms, natural in the second argument (that is coherent with
morphisms of these second arguments). Here and thereafter ${\bf h}_{(\cd)}(\cd,\cd)$ denotes the set
of all morphisms between respective objects of a category in question.

Suppose that a given rig has such a nice property: every object in ${\cal L}$ has a free object in 
${\cal K}$ with that base. In this case our rigged category is sometimes called {\it 
freedom-loving}. When it happens, we can, for every object $X$ in ${\cal K}$, apply the map 
 ${\cal I}_{\sq(X),X}$ to the identity morphism $\id_{\sq(X)}\in{\bf h}_{{\cal L}}(\sq(X),\sq(X))$.
 The resulting morphism $\pi_X:{\bf Fr}(\sq(X))\to X$ is called the {\it canonical morphism for $X$}.
 Then a categorical argument, actually, contained in~\cite{mc2}, leads to the following statement.
 We shall use it essentially in the next section.

\medskip
{\bf Proposition 3}. {\it In the case of a freedom-loving rigged category an object $P$ in ${\cal 
K}$ is $\sq$-projective if and only if the canonical morphism $\pi_P$ has a right inverse morphism 
in ${\cal K}$.} 

\medskip

As a matter of fact, all rigged categories, providing three above-mentioned types of the
projectivity, are freedom-loving.
 However, we restrict ourselves with the rig (6), providing the relative projectivity.

Our main tool is the notion of the so-called operator-projective tensor product of Q-spaces, independently
discovered in~\cite{er2} and~\cite{blp}. Following~\cite{heb2}, we shall define it in terms of  its
universal property.

Let us fix, for a time, Banach Q-spaces $E$ and $F$.

\medskip
{\bf Definition 8.}  A pair $(\Theta,\theta)$, consisting of a Banach Q-space $\Theta$ and a
completely contractive bilinear operator $\theta:E\times F\to\Theta$, is called (completed) {\it
Q-tensor product of $E$ and $F$} if, for every completely bounded bilinear operator
$\rr:E\times F\to G$, where $G$ is a Banach Q-space, there exists a unique completely bounded
operator $R:\Theta\to G$ such that the diagram
\[
\xymatrix@R-10pt@C+15pt{
E\times F \ar[d]^{\theta} \ar[dr]^{\rr} & \\
\Theta \ar[r]^R  &  G  }
\]
\noindent is commutative, and, moreover, $\|R\|_{cb}=\|\rr\|_{cb}$.

\medskip
Such a pair does indeed exist. Here we only recall, without proof, its explicit construction in the frame-work of the non-coordinate presentation. As we shall see, $\Theta$
turns out to be the completion of a certain Q-space $E\ot_{op}F$, which is $E\ot F$ with respect to a special Q-norm; we shall denote this completion by $E\widehat\ot_{op}F$. As to $\theta$, it is
 just the canonical bilinear operator $\vartheta:E\times F\to E\ot F:(x,y)\mt x\ot y$, only
considered with the range space $E\widehat\ot_{op} F$.

To introduce the mentioned Q-norm, we need some ``extended'' version of the diamond multiplication,
this time between elements of the amplifications of linear spaces. Namely, for $u\in \cK E,v\in\cK F$
we denote by $u\di v$ the element $\vartheta_\ii(u,v)\in\cK(E\ot F)$. In other words, this kind of ``diamond
operation'' is well defined on elementary tensors by $ax\di by:=(a\di b)(x\ot y).$

The first observation is that {every $U\in\cK(E\ot F)$ can be represented as
$$
a\cd(u\di v)\cd b\eqno (8)
$$
\noindent for some $a,b\in\cK, u\in\cK E,v\in\cK F$. For a simple proof see, e.g.,~\cite[Prop. 7.2.10]{heb2}.

\medskip
After this, for every $U$ we introduce the number
$$
\|U\|_{op}:=\inf\left\{\|a\|\|u\|\|v\|\|b\|\right\}, \eqno (9)
$$
where the infimum is taken over all possible representations of $U$ in the form (8). The following
theorem is proved as Theorem 7.2.19 in~\cite{heb2}.

\medskip
{\bf Theorem.} {\it The function $U\mt\|U\|_{op}$ is a Q-norm on $E\ot F$. Further, if we denote the respective Q-space by
$E\ot_{op}F$ and its completion by $E\widehat\ot_{op}F$, then the pair $(E\widehat\ot_{op}F,\vartheta)$
is a (completed) Q-tensor product of Banach Q-spaces $E$ and $F$.}

\medskip
As a part of this assertion, $\vartheta$ is completely contractive, that is
$$
\|u\di v\|_{op}\le\|u\|\|v\|; \q u\in \cK E,v\in\cK F.\eqno (10)
$$

Further,
by (3) we have that $\|a\di b\|=1$ provided $\|a\|=\|b\|=1$. Therefore the
action of `` $\di$ ' on elementary tensors (see above), combined with (1),
implies that in the normed space $E\widehat\ot_{op}E$ we have
$$
\|x\ot y\|\le\|x\|\|y\|; \q x\in E, y\in F. \eqno  (11)
$$
 (In fact,in (10) and (11) we have the exact equality, but we shall not discuss it now).

\medskip
What does this tensor product give for the questions concerning projectivity and freeness?

In what follows, the dot `` $\cd$ ''
always denotes the outer multiplication, whatever base algebra is considered at the moment. This
will not create a confusion. Also we shall denote the norm on $E\widehat\ot_{op}E$ just by $\|\cd\|$.

\medskip
{\bf Proposition 4}. {\it Let $A$ be a Banach Q-algebra, $X$ a left Banach Q-module over $A$,
$E$ a Banach Q-space. Then the Banach Q-space $X\widehat\ot_{op} E$ has a unique structure of a left
Banach Q-module over $A$, such that for elementary tensors in $X\ot E$ we have
${b}\cd(x\ot y)=({b}\cd x)\ot y$; ${b}\in A, x\in X, y\in E$.}

\smallskip
 The proof could be deduced more or less quickly from the associativity property of the operation of the
operator-projective tensor product. However, the complete proof of the said property
is rather long and technical (of course, if it is not prudently left to the reader,
as in~\cite[Prop. 7.1.4]{efr}).
Therefore we prefer to give an independent proof.

$\tl$ The uniqueness of the indicated structure immediately follows from the density of $X\ot_{op}E$ in
$X\widehat\ot_{op} E$. We proceed to the proof of its existence.

Fix, for a time, ${\bf a}\in A$ and set
${\cal S}:X\times E\to X\widehat\ot_{op} E:(x,y)\mt({\bf a}\cd x)\ot y$.
Since  the bilinear operator of the outer multiplication in $X$, denoted by $\widetilde{\cal M}$, is completely bounded, the ``partial'' operator $^{\bf a}\widetilde{\cal M}:X\to X:x\to{\bf a}\cd x$
is also completely bounded, and $\|^{\bf a}\widetilde{\cal M}_\ii\|\le\|{\bf a}\|$;
this is by virtue of Proposition 1. Further, we have the formula
$$
{\cal S}_\ii(u,v)=\vartheta_\ii(^{\bf a}\widetilde{\cal M}_\ii u,v); u\in\cK E, v\in\cK F.
$$
(LOOK AT ELEMENTARY TENSORS). Therefore we have $\|{\cal S}_\ii(u,v)\|\le\|^{\bf a}\widetilde{\cal M}_\ii u\|\|v\|\le\|{\bf a}\|\|u\|\|v\|$. Thus ${\cal S}$ is completely bounded with $\|{\cal S}_{cb}\|\le\|{\bf a}\|$, and hence gives rise to the completely bounded operator $^{\bf a} \frak {M}:X\widehat\ot_{op} E\to X\widehat\ot_{op} E$ with the same estimation of its completely bounded norm. Obviously we have $^{\bf a} \frak   {M}(x\ot y)=(a\cd x)\ot y$.

Now ``release'' ${\bf a}$ and set $\frak {M}:A\times(X\widehat\ot_{op} E)\to X\widehat\ot_{op} E:({\bf a},z)\mt^{\bf a}\frak{M}(z);\\ z\in X\widehat\ot_{op} E$. Obviously, $\frak {M}$ is a left outer multiplication in $X\widehat\ot_{op} E$, acting on elementary tensors as it is indicated in the formulation. Our task is to show that it is completely contractive, that is for every $w\in\cK A$ and $V\in\cK X\widehat\ot_{op} E$ we have $\|\frak {M}_\ii(w,V)\|\le\|w\|\|V\|$.

First, it is easy to observe that there is a unitary operator $T$ on $L$ such that for every $a,b,c\in\bb$ we have $(a\di b)\di c=T[a\di(b\di c)]T^*$. (see~\cite[p. 15]{heb2} for details of the proof). This  implies the formula
$$
\vartheta_\ii(\widetilde{\cal M}_\ii(w,u),v)=T\cd [\frak {M}_{\ii}(w,\vartheta_{\ii}(u,v))]\cd T^*;w\in\cK A, u\in\cK X,v\in\cK E.
$$
(LOOK AT ELEMENTARY TENSORS). But $\widetilde{\cal M}$ and $\vartheta$ are completely contractive, and it follows from the first Ruan axiom for Q-spaces that $\|T\cd V\cd T^*\|=\|V\|$ for all $V\in X\widehat\ot_{op} E$. Therefore we have
$$
\|\frak {M}_\ii(w,u\di v)\|\le \|w\|\|u\|\|v\|.
$$
Further, take an arbitrary $U\in\cK(X\ot E)$ and represent it as $a\cd(u\di v)\cd b$ (cf. (8)). Observe
the formula
$$
\frak {M}_\ii(w,U)=(\id\di a)\cd\frak {M}_\ii(w,u\di v)\cd(\id\di b).
$$
(LOOK AT ELEMENTARY TENSORS). This, combined with (4), implies that
$$
\|\frak {M}_\ii(w,U)\|\le\|\id\di a\|\|\frak {M}_\ii(w,u\di v)\|\|\id\di b\|=
$$
$$
\|a\|\|\frak {M}_\ii(w,u\di v)\|\|b\|\le\|w\|\|a\|\|u\|\|v\|\|b\|.
$$
Taking all representations of $U$ in the indicated form and recalling the definition of the norm on
$\cK(X\ot_{op}E)$ (see (9)), we get the estimation $\frak {M}_\ii(w,U)\|\le\|w\|\|U\|$.

It remains to show that such an estimation holds for all $V\in\cK(X\widehat\ot_{op}E)$, not necessary belonging to $\cK(X\ot E)$. Since the latter space is dense in $\cK(X\widehat\ot_{op}E)$, it is sufficient to show that for every $w\in\cK A$ the function $V\mt\|\frak {M}_\ii(w,V)\|$ is continuous.
But elements in $\cK A$ are sums of elementary tensors, and we have the triangle inequality for the norm. 
Therefore it suffices to show that for all $b\in\cK,{\bf a}\in A$
the element $\frak {M}_\ii(b{\bf a},V)\in\cK(X\widehat\ot_{op}E)$ continuously depends on $V$.

To prove the latter assertion, we notice that for an elementary tensor in $\cK(X\widehat\ot_{op}E)$, say
$c\Phi; c\in\cK,\Phi\in X\widehat\ot_{op}E$ we have
$$
\frak {M}_\ii(b{\bf a},c\Phi)=(b\di c)\frak {M}({\bf a},\Phi)=(b\di c)[^{\bf a}\widetilde{\cal M}(\Phi)]=
b\di (c[^{\bf a}\widetilde{\cal M}(\Phi)])=b\di[^{\bf a}\widetilde{\cal M}_\ii(c\Phi)].
$$
It follows, by the additivity of relevant operations, that we have
$$
\frak {M}_\ii(b{\bf a},V)=b\di[\widetilde{\cal M}^{\bf a}_\ii(V)]\eqno (12)
$$
for all $V\in\cK(X\widehat\ot_{op}E)$. But we remember that the operator $^{\bf a}\widetilde{\cal M}$,
acting on $X\widehat\ot_{op}E$, is completely bounded. This obviously implies that
$^{\bf a}\widetilde{\cal M}_\ii(V)\in\cK(X\widehat\ot_{op}E)$ continuously depends on $V$.
Consequently, taking into account (12) and (4), we see that the same is true for
$\frak {M}_\ii(b{\bf a},V)$. And this, as it was said above, is just what we need. $\tr$

\medskip
As the most important particular case of the latter proposition, we can take, in the place of $X$,
the base algebra $A$ and speak about the left Q-module $A\widehat\ot_{op}E$ with the outer multiplication,
well defined by ${a}\cd(b\ot y)=(ab)\ot y$.

\medskip
From now on, for simplicity, we suppose that a given Banach Q-algebra $A$ has an identity of norm 
1, denoted in what follows by ${\bf e}$, and that we have the identity ${\bf e}\cd x=x$ for all 
modules in question. 

\medskip
Otherwise we would deal with the unitization of $A$, which itself can be made into a Banach 
Q-algebra with respect to a certain Q-norm on $A\oplus\co$. The latter is a particular case of the 
so-called operator, or quantum $l_1$-sum; see~\cite{ble,pis}. But we shall not give the details. 

\bigskip
This is the reason why we need such a module:

\medskip
{\bf Proposition 5}. {\it Under the given assumptions about $A$ and $E$, the Q-module
$A\widehat\ot_{op} E$, being considered  in the rigged category $({\bf QA\mbox{-}mod},\sq)$,
 is a free object with the base $E$.  Namely,} (according to the general definition of the freeness,
 given above; cf. (7)) {\it  for every Q-module $Y$ over $A$ there is a bijection
$$
{\cal I}_{E,Y}:{\bf h}_{{\bf QBan}}(E,\square Y)\to{\bf h}_{{\bf QA-mod}}(A\widehat\ot_{op}E,Y),
$$
natural in $Y$. This bijection takes  a completely bounded operator $\va:E\to Y$ to the morphism
$S:A\widehat\ot_{op}E\to Y$, well defined by ${\bf a}\ot x\mt{\bf a}\cd\va(x)$ and natural in $Y$.}

\smallskip
$\tl$ Let $\va$ be as before. Consider the bilinear operator ${\cal S}:A\times E\to Y$, taking a pair
$({\bf a},x)$ to ${\bf a}\cd\va(x)$, or, what is the same, to $\widetilde{\cal M}({\bf a},\va(x))$,
where $\widetilde{\cal M}:A\times Y\to Y$ denotes the respective bilinear operator of outer multiplication.
 Observe the formula
$$
{\cal S}_\ii(w,u)=\widetilde{\cal M}_\ii(w,\va_\ii(u)); w\in\cK A, u\in\cK E.  \eqno (13)
$$
(LOOK AT ELEMENTARY TENSORS). It follows that $\|{\cal S}_\ii(w,u)\|\le\|w\|\|\va_\ii(u)\|\le\|\va\|_{cb}\|w\|\|u\|$. Therefore
${\cal S}$ is completely bounded, and $\|{\cal S}\|_{cb}\le\|\va\|_{cb}$. Consequently, ${\cal S}$ gives
rise to the completely bounded operator $S:A\widehat\ot_{op}E\to Y$, well defined as it is indicated in
the formulation.

 It follows that
$$
S({\bf a}\cd\Phi)={\bf a}\cd S(\Phi)\eqno (14)
$$
for all ${\bf a}$ in $A$ and all $\Phi$ in the dense subspace $A\ot E$ of $A\widehat\ot_{op}E$.
( LOOK AT ELEMENTARY TENSORS). But both parts of this equality continuously depend on $\Phi$ when the latter
runs the whole $A\widehat\ot_{op}E$: this immediately follows already from the ``classical'' boundedness of
$\widetilde{\cal M}$ and $S$
Consequently, we have an equality (13) for all $\Phi\in A\widehat\ot_{op}E$; in other words, $S$ is a
morphism of $A$-modules.

Thus we obtain a map ${\cal I}_{E,Y}:\va\mt S$ between the sets, indicated in the formulation.
 A similar argument, again using the density of $A\ot E$ of $A\widehat\ot_{op}E$ and the boundedness of
 relevant bilinear and linear operators, shows that this map is natural in $Y$.

It remains to show that ${\cal I}_{E,Y}$ is a bijection. For this aim we shall display its inverse map. Take an arbitrary morphism $S:A\widehat\ot_{op}E\to Y$ and consider the operator $\va:E\to Y:x\mt S({\bf e}\ot x)$.
Of course, $\va$ is not other thing that the composition $S[^{\bf e}\vartheta]$, where
$^{\bf e}\vartheta:E\to A\widehat\ot_{op}E$ is the relevant ``partial'' operator
with respect to $\vartheta$
 Since $S$ and, by Proposition 1, $^{\bf e}\vartheta$ are completely bounded, the same is obviously
 true for $\va$.  Assigning such a $\va$ to every  $S$, we obtain a map
 $J_{E,Y}:{\bf h}_{{\bf QA-mod}}(A\widehat\ot_{op}E),Y)\to{\bf h}_{{\bf QBan}}(E,\square Y)$. From the
 definitions of $I_{E,Y}$ and $J_{E,Y}$ one can easily see that the compositions
 ${\cal j}_{E,Y}{\cal I}_{E,Y}$ and ${\cal I}_{E,Y}{\cal J}_{E,Y}$ are identity maps on the respective
 sets of morphisms. This completes the proof. $\tr$

\medskip
Now recall the notion of the canonical morphism, defined above in the frame-work of a general
freedom-loving rigged category. How does it act in the case of our special rig (7)?
 Take a Q-module $X$ over $A$. We know the special form of the bijection ${\cal I}_{E,Y}$, indicated
 in the previous proposition. Setting $E:=\sq(X),Y:=X$, we immediately see that the canonical morphism
 for $X$ is $\pi_X:A\widehat\ot_{op}X\to X$, well defined on elementary tensors by taking ${\bf a}\ot x$ to
${\bf a}\cd x$. (Here, of course, in the module $A\widehat\ot_{op}X$ we consider $X$ just as a 
Q-space). 

Consequently, as a particular case of Proposition 3, we obtain

\medskip
{\bf Proposition 6}. {\it A Q-module $P$ over a Q-algebra $A$ is relatively projective if and only 
if the canonical morphism $\pi_P$ has a right inverse (completely bounded) $A$-module morphism.} 

\bigskip
{\bf 3. Quantum projectivity and ideals in $C^*$-algebras}

\bigskip
The main result of this section is

\medskip
{\bf Theorem.} {\it Let $A$ be a separable $C^*$-algebra, endowed with the standard quantization} 
(see Example), {\it $I$ be a closed left ideal of $A$. Then $I$, considered as left Q-module over 
$A$, is relatively projective.} 

 \smallskip
As it was said in Introduction, the ``classical'' prototype of this theorem was obtained by Z.A.Lykova.

\smallskip
Roughly speaking, our argument consists of two parts: ``classical'' and ``quantum''. As to the
first part, it resembles what was done in~\cite{he1,lyk,phr}. However,
we shall use somewhat sharper estimation of norms of certain elements of the ideal in question.

\smallskip
$\tl\tl$
It is well known (cf., e.g., Sections 1.7.2,1.7.3 in~\cite{dix}) that
 $I$ has a positive countable left approximate identity of norm $<1$, denoted in what follows by
 $e_n; n\in{\B N}$. Taking the $C^*$-algebra,
 generated  by elements $e_n$, and applying to it Corollary 1.5.11 in~\cite{lin}, we can assume that
 $e_n$, in addition, is such that
 $e_ne_{n+1}=e_{n+1}e_{n}=e_{n}$. Set $e_0:=0$ and $b_n:=\sqrt{e_{n}-e_{n-1}}; n\in{\B N}$.

\medskip
In the following lemma we are given $x\in I, m,n\in{\B N}; m<n$, and also, for $k=m+1,...,n$, complex numbers $\xi_k; |\xi_k|=1$. Set $C:=\max\{\|x-xe_k\|;k=m,...,n+1\}$.

\medskip
{\bf Lemma.} {\it  We have

\medskip
(i) $\|\sum_{k=m+1}^nx\xi_kb_k\|\le\sqrt{6C\|x\|}$.

\medskip
(ii) $\|\sum_{k=m+1}^n\xi_kb_k\|\le\sqrt{2}$. }

\medskip
$\tl$ \' (i) Denote the first sum by $z$. Because of $C^*$-property, we have $\|z\|^2=\|zz^*\|$.

\smallskip
We see that $zz^*=B_1+B_2$, where $B_1:=\sum_{k=m}^nxb_k^2x^*=x(e_{n}-e_{m})x^*$, and

 \smallskip
 $B_2:=\sum\{x[\xi_k\bar\xi_l b_kb_l]x^*; k,l=m+1,...n, k\ne l\}$ in the case
$m+1<n$ whereas $B_2:=0$ otherwise. An easy calculation shows that
$$
\|B_1\|=\|[(xe_n-x)-(xe_m-x)]x^*\|\le[\|xe_{n}-x\|+\|x-xe_{m}\|]\|x^*\|\le2C\|x\|.\eqno (15)
$$

\medskip
Turn to $B_2$. It follows from the choice of $b_n$ and properties of $e_n$ that
$$
b_kb_l=0 \q\q {\rm whenever} \q\q |k-l|>1.\eqno (16)
$$
Therefore we have
$$
B_2=\sum_{k=m+1}^{n-1}x\xi_kb_kb_{k+1}\bar\xi_{k+1}x^*+
\sum_{k=m+1}^{n-1}x_{k+1}\xi_{k+1}b_{k+1}b_{k}\bar\xi_kx^*=
\sum_{k=m+1}^{n-1}xt_kb_kb_{k+1}x^*,
$$
where
$t_k:=\xi_k\bar\xi_{k+1}+\xi_{k+1}\bar\xi_{k}$. In particular, we see that $B_2$, as well as, of
course, $B_1$, is self-adjoint.

Now note that we have $-2\le t_k\le2$, and also, with respect to the order in $A$, the estimate
$2b_kb_{k+1}\le b_k^2+b_{k+1}^2=e_{k+1}-e_{k-1}$ holds. Consequently, for
all $k=m+1,...,n-1$ we have
$$
-(e_{k+1}-e_{k-1})\le-2b_kb_{k+1}\le t_kb_kb_{k+1}\le2b_kb_{k+1}\le e_{k+1}-e_{k-1}.\eqno(17)
$$
Summing these inequalities and multiplying the resulting inequality from different sides by
$x$ and $x^*$, we obtain that
$$
-x[e_{n-1}+e_{n}-e_{m-1}-e_{m}]x^*\le B_2\le x[e_{n-1}+e_{n}-e_{m-1}-e_{m}]x^*,
$$
and hence
$$
\|B_2\|\le\|x[e_{n-1}+e_{n}-e_{m-1}-e_{m}]x^*\|=
$$
$$
\|[(xe_{n-1}-x)+(xe_{n}-x)-(xe_{m-1}-x)-(xe_{m}-x)]x^*\|\le4C\|x\|
$$

\smallskip
Hence $\|z\|^2\le\|B_1\|+\|B_2\|\le6C\|x\|$, and we are done.

\medskip
(ii) A similar estimate was already obtain in~\cite{phr} with the help of the commutative
Gelfand/Naimark Theorem. Instead, we choose to prove it, using some properties of norms in $C^*$-algebras.

Again, denote the relevant sum by $z$. Now $zz^*=e_{n}-e_m+B$, where in the case $m+1<n$ we have,
because of (16),
$B=\sum_{k=m+1}^{n-1}t_kb_kb_{k+1}$, and otherwise we have $B=0$; here $t_k$ is the
same as in (i). Further, we obtain from (17) that $\|t_kb_kb_{k+1}\|\le\|e_{k+1}-e_{k-1}\|=1$
for all $k$. Finally, it immediately follows from (16) that $b_kb_{k+1}b_lb_{l+1}=0$, provided
$k\ne l$. Therefore for every natural $N$ we have either
$$
\|B\|=\|B^{2^N}\|^{\frac{1}{2^N}}=
\left\|\left(\sum_{k=m+1}^{n-1}(t_kb_kb_{k+1})^{2^N}\right)\right\|^{\frac{1}{2^N}}\le
$$
$$
\left(\sum_{k=m+1}^{n-1}\|t_kb_kb_{k+1}\|^{2^N}\right)^{\frac{1}{2^N}}\le(n-m-1)^{\frac{1}{2^N}}.
$$
or $B=0$. Hence $\|B\|\le1$. Therefore $\|z\|^2\le\|e_{n}-e_m\|+1\le2$. $\tr$

\medskip
{\it Continuation of the proof: new estimates}

\medskip
Retain the notation of the previous lemma. Since $b_k\in I$, we can consider in the free Q-module $A\widehat\ot_{op}I$
elements $V_{m,n}:=\sum_{k=m}^n xb_{k}\ot b_{k}$ and $W_{m,n}:=\sum_{k=m}^n b_{k}\ot b_{k}$.

\smallskip
Recall an old trick, used in  arithmetical Lemma 2.41 in~\cite{he1}. Namely, consider the complex number $\zeta:=e^{2\pi i/(n-m+1)}$, that is the $(n-m+1)$-th primitive root of unity. A routine calculation (cf. {\it idem}) shows that we have
$$
V_{m,n}=\frac{1}{n-m}\sum_{k=1}^{n-m}\left[\sum_{i=m+1}^nx\zeta^{k(i-n-1)}b_i\right]\ot
\left[\sum_{j=m+1}^n\zeta^{-k(j-n-1)}b_j\right]
$$
and
$$
W_{m,n}=\frac{1}{n-m}\sum_{k=1}^{n-m}\left[\sum_{i=m}^n\zeta^{k(i-n-1)}b_i\right]\ot
\left[\sum_{j=m}^n\zeta^{-k(j-n-1)}b_j\right].
$$
The estimate (11), combined
with the previous lemma, implies that
$$
\|V_{m,n}\|\le2\sqrt{3C\|x\|} \q\q {\rm whereas} \q\q \|W_{m,n}\|\le2.\eqno (18)
$$

\medskip
{\it Continuation of the proof: The appearance of the morphism $\rho$}

\medskip
Now, remembering Proposition 6, we proceed to the construction a completely bounded morphism
$\rho:I\to A\widehat\ot_{op} I$, right inverse to $\pi_I$.

Consider, for every $x\in I$ and $n\in{\B N}$, the element
$$
\rho_n(x):=\sum_{k=1}^n xb_{k}\ot b_{k}\in A\widehat\ot_{op} I,
$$
and the resulting sequence of operators $\rho_n:I\to A\widehat\ot_{op} I:x\mt\rho_n(x)$. Fix, for a moment, $n$ and set, for brevity,
$W:=W_{1,n}$ (see above). Then, in the notation of Proposition 1, we obviously have
$\rho_n(x)=\widetilde{\cal M}^{W}(x)$, where $\widetilde{\cal M}:A\times I\to I$ is the bilinear operator
of the outer multiplication in the $A$-module $I$. Therefore, by the mentioned proposition, $\rho_n$ is
completely bounded, and $\|\rho_n\|_{cb}\le\|W\|$. Hence, by (18), we have $\|\rho_n\|_{cb}\le2$.

 Fix $\e>0$. Then there exists natural $N$ such that $\|x-xe_n\|<{\e}/(2\sqrt{3}\|x\|)$
 whenever $n>N$. Take $m,n; N<m<n$. Then $\rho_n(x)-\rho_m(x)=\sum_{k=m+1}^n xb_{k}\ot b_{k}$. Therefore,
 by (17), $\|\rho_n(x)-\rho_m(x)\|<\sqrt{3}\e\|x\|/(2\sqrt{3}\|x\|)=\e$. Thus we see that $\rho_n(x)$ is a
 Cauchy net, and hence it converges in the Banach space $A\widehat\ot_{op}I$ to some element; denote the
 latter by $\rho(x)$. In this way the map $\rho:I\to A\widehat\ot_{op}I:x\mt\rho(x)$ appears. It is
 easy to see that $\rho$ is a morphism of $A$-modules.

To move further, let us distinguish a statement of a general character.

\medskip
{\bf Lemma.} {\it Suppose that $E,F$ are Q-spaces, and $T_n:E\to F, n\in{\B N}$ are completely bounded
operators. Suppose also that for every $x\in E$ the sequence $T_n(x)$ converges to some $Tx\in F$, and
there is $C>0$ such that $\|T_n\|_{cb}<C$ for all $n$. Then the map $T:E\to F:x\mt Tx$ is also a
completely bounded operator, and $\|T\|_{cb}\le C$.}

\smallskip
$\tl$ Of course, $T$ is an operator. Take $u\in{\cal F}E; u=\sum_{k=1}^na_kx_k$. Then
$(T_n)_\ii(u)=\sum_{k=1}^na_kT(x_k)$, and therefore
the equality (1)
implies that $(T_{n})_\ii(u)$ converges to $T_\ii(u)$.
Consequently the estimate $\|(T_{n})_\ii(u)\|\le
C\|u\|$ implies that $\|T_\ii(u)\|\le C\|u\|$, and we are done. $\tr$

\medskip
{\it The end of the proof}

\medskip
Since $\|\rho_n\|_{cb}\le2$, the previous lemma implies that $\rho$
is also completely bounded with the same estimate $\|\rho\|_{cb}\le2$.

It remains to show that $\pi_I\rho=\id_I$. Indeed, for every $x\in I$ we have
$$
\pi_I\rho(x)=\pi_I(\lim_{n\to\ii}\rho_n(x))=
\lim_{n\to\ii}\pi_I\left(\sum_{k=1}^n xb_{k}\ot b_{k}\right)=
 \lim_{n\to\ii}\left(\sum_{k=1}^nxb_k^2\right)=
 $$
 $$
\lim_{n\to\ii}\left[\sum_{k=1}^n(xe_{k}-xe_{k-1})]\right]=
 \lim_{n\to\ii}xe_{n}=x. \q\q
$$
The theorem is proved. $\tr\tr$

\medskip
{\bf Remark.} The condition of separability of $A$ can be weakened. In particular, the result is valid
provided our algebra has a strictly positive element, that is $a\ge0$ with $\om(a)>0$ for all states
$\om$ on $A$. The proof is practically the same as given above.

On the other hand, ideals in general $C^*$-algebras are not bound to be relatively projective.
For example, suppose that our algebra is commutative, and the Gelfand spectrum of a given
ideal is not paracompact, like in the case of maximal ideals in $C(\beta{\B N})$, corresponding to points
of $\beta{\B N}\setminus{\B N}$. Then
 in the ``classical'' context such an ideal is not projective (cf. what was said in Introduction).
 The same argument, up to minor modifications, shows that the same is true in the ``quantum'' context.

\medskip
{\bf Remark.} We have already mentioned that passing from the relative to topological and metric
projectivity we get much less projective modules. In this connection we would like to cite
a rather difficult theorem, due to
N.Nemesh~\cite{nor} and concerning the ``classical'' context. Namely, Nemesh proved
that for a closed left ideal, say $I$, in a $C^*$-algebra the following properties are equivalent:

(i) $I$ is topologically projective

(ii) $I$ is metrically projective

(iii) $I$ has a right identity which is a self-adjoint idempotent

\medskip
We believe that such a theorem holds in the ``quantum'' context as well, although so far we have not seen
an accurate proof.

\ed

{\bf Theorem 4.} {\it Let $I$ be a closed left ideal in a $C^*$-algebra $A$. Then its following
properties are equivalent:

(i) $I$ is topologically projective

(ii) $I$ is metrically projective

(iii) $I$ has a right identity which is a self-adjoint idempotent

(iv) $I$ has a complementing closed left ideal, say $J$, in $A$ such that the projection of $A$
onto $I$ along $J$ is contractive.}

\medskip
Let us stress that in both cases $I$ necessarily has a right identity or, equivalently, has a
closed module complement in $A$.

This theorem was formulated in the ``classical'' context. Most probably, it has a natural version
for quantum modules, but still there is no accurate proof.

\bigskip
In conclusion, I shall briefly report some results, concerning two other types of ``homologically
best'' modules, injectivity and flatness. Both are formulated in the classical context, where the
definitions of the injectivity and flatness are well known for my present audience.

The first theorem concerns the most important class of von Neumann algebras that has many names:
Connes-amenable, hyperfinite, semi-discrete etc. An old result of 1989 states that a von Neumann
algebra $A$ is hyperfinite if, and only if its predual $A$-bimodule $A_*$ is relatively injective.
In 2010 G.Racher~\cite{rah} obtained a more precise ``quantitative'', in the spirit of the paper of
M.White~\cite{whi} of 1996, form of this result:

\medskip
{\bf Theorem 5.} {\it A  von Neumann algebra $A$ is hyperfinite iff its predual $A$-bimodule $A_*$
is relatively injective (that is for every extension problem for its predual bimodule such that the
relevant natural embedding has a right inverse $j$, there is its solution with the norm
$\le\|\va\|\|j\|$).}

\bigskip
Finally, let us proceed to the third fundamental homological notion, the flatness. You know the
definition; I only recall that it is more tolerant and widespread, and sometimes more practical
concept that the projectivity.

Let us concentrate on cyclic modules, that is left $A$-modules of the form $X:=A/I$, where $A$ is a
Banach algebra, now for simplicity supposed to be unital, and $I$ its closed left ideal. The old
result of 1979, due to M.Sheinberg and the speaker, claims that for the relative flatness of such a
module it is sufficient that $I$ have a right bounded approximate identity, and if $I$ is weakly
complemented in $A$, the converse is also true. (The latter condition from the Banach geometry can
not be omitted, and this, curiously enough, led to the construction of new examples of
non-complemented subspaces of Banach spaces.) However, if we pass from the relative to the
topological flatness, we obtain an exact criterion. The following result modifies a certain
``quantitative'' assertion of White (see {\it idem}), and in the following form is presented by
 Nemesh.

 \medskip
 {\bf Theorem 6.} {\it A cyclic module $A/I$
$X$ is topologically flat if, and only if $I$ has a right bounded approximate identity.}